\documentclass[12pt,twoside]{article}

\date{25 July 2008}

\newtheorem{Theorem}{\quad Theorem}[section]

\newtheorem{Definition}[Theorem]{\quad Definition}

\newtheorem{Lemma}[Theorem]{\quad Lemma}

\newcommand{\ds}{\displaystyle}

\begin{document}

\centerline{\bf Uzbek Mathematical Journal,  1996, No.4 }

\centerline{}

\centerline{}

\centerline {\Large{\bf About the asymptotic formula for spectral
}}

\centerline{ }

\centerline{\Large{\bf function of the Laplace-Beltrami operator
on sphere }}

\centerline{}

\newcommand{\mvec}[1]{\mbox{\bfseries\itshape #1}}

\centerline{\bf {Anvarjon Akhmedov}}

\centerline{Tashkent State University,}

\centerline{Department of Mathematical Prysics,}

\centerline{ Vuzgorodok, Tashkent, Uzbekistan.}

\centerline{ anv.akhmed@gmail.com}

\begin{abstract}
\textbf{\emph{In this work we established asymptotical behavior
for Riesz means of the spectral function of the Laplace operator
on unit sphere.}}
\end{abstract}


\section{Introduction}

Let $S^N$ is the $N$-dimensional unit sphere in $R^{N+1}$. The
Laplace operator $\Delta$:
$$
\Delta u=\frac{\partial^2u}{\partial x_1^2}+
\frac{\partial^2u}{\partial x_2^2}+...+\frac{\partial^2u}{\partial
x_N^2}, u\in C_0^\infty(S^N)$$

in spherical coordinates in $N$ dimensions, with the
parametrization $x=r\theta \in  R^{N+1}$ with $r \in [0,+\infty)$
and $\theta \in S^{N}$,
$$
  \Delta f = \frac{\partial^2 f}{\partial r^2} + \frac{N}{r} \frac{\partial f}{\partial r}
  + \frac{1}{r^2} \Delta_Sf
$$
where $\Delta_{S}$ is the Laplace-Beltrami operator on the $N$-
 dimensional sphere, or spherical Laplacian. One can also write
the term
$${\partial^2 f \over \partial r^2} + \frac{N}{r}
\frac{\partial f}{\partial r}
$$

equivalently as

$$\frac{1}{r^{N}} \frac{\partial}{\partial r}
\Bigl(r^{N} \frac{\partial f}{\partial r} \Bigr).
$$

As a consequence, the spherical Laplacian of a function defined on
$S^{N}\subset{ R}^{N+1}$ can be computed as the ordinary Laplacian
of the function extended to ${R}^{N+1} \setminus\{0\}$ so that it
is constant along rays.

Let us denote by $\lambda_0, \lambda_1,...$ the distinct
eigenvalues of $-\Delta_S$, arranged in increasing order. Let
$H_k$ denote  the eigenspace corresponding to $\lambda_k.$ We call
elements of $H_k$ spherical harmonics of degree $k$. It is well
known (see \cite{ST5}) that $dimH_k=a_k:$
\begin{eqnarray}\label{eq12*}
a_k=\left\{
\begin{array}{ll}

\ds \qquad 1,\ \qquad\qquad\qquad\qquad\qquad if \qquad k=0,\\[2mm]

\ds \qquad N , \qquad \qquad\qquad\qquad\qquad if \qquad k=1, \\[2mm]

\ds \frac{(N+k)!}{N!k!}-\frac{(N+k-2)!}{N!(k-2)!},  \qquad if
\qquad k\geq 2 \\

\end{array}
\right.
\end{eqnarray}

For $a_k$ we have $a_k\approx k^{N-1}$ as $k\rightarrow\infty.$

Let $\hat{A}$ is a self-adjoint extension of the Laplace operator
$\Delta_{S}$  in $L_2(S^N)$ and if $E_\lambda$ is the
corresponding spectral resolution, then for all functions $f\in
L_2(S^N)$ we have

$$
\hat{A}f=\int_0^{\infty}\lambda dE_\lambda f.
$$

The operator$\hat{A}$ has in$L_2(S^N)$ a complete orthonormal
 system of eigenfunctions
  $$ \{Y_1^{(k)}(x),Y_2^{(k)}(x),...,Y_{a_k}^{(k)}(x)\}\subset H_k, k=0,1,2,...,
$$
corresponding to the eigenvalues $\{\l_k=k(k+N-1)\}, k=0,1,2,...$.

 It is easy to check
that the operators $E_\lambda$ have the form
\begin{equation}\label{PS}
E_\lambda f(x)=\sum_{\lambda_n<\lambda}Y_n(f,x) ,
\end{equation}
 where
\begin{equation}\label{FC}
 Y_k(f,x)=\sum\limits_{j=1}^{a_k}Y_j^{(k)}(x)\int\limits_{S^N}f(y)Y_j^{(k)}(y)d\sigma(y)
 \end{equation}.

  The Riesz means of the partial sums (\ref{PS}) is defined by

\begin{equation}\label{RM}
E_n^\alpha
f(x)=\sum\limits_{k=0}^n\left(1-\frac{\lambda_k}{\lambda_n}\right)^\alpha
Y_k(f,x)
\end{equation}

The most convenient object for a detailed investigation are the
expansions of the form  (\ref{RM}). The  integral (\ref{RM}) may
be transformed writing instead of $\hat{A}$ the integral to the
right in (\ref{FC}) and then changing the order of integration.
This yields the formula
 \begin{equation} \label{RM1}
 E_\lambda^s
 f(x)=\int\limits_{S^N}\Theta^s(x,y,\lambda)f(y)d\sigma(y)
 \end{equation}
with

\begin{equation}\label{05}
\Theta^\alpha(x,y,n)=\sum\limits_{k=0}^n\left(1-\frac{\lambda_k}{\lambda_n}\right)^\alpha
Z_k(x,y).
\end{equation}

For $\alpha=0$ this kernel is called the spectral function of the
Laplace operator for the entire space $S^N$.

The behavior of the spectral expansion corresponding to the the
Laplace-Beltrami operator is closely connected with the
asymptotical behavior of the kernel $\Theta^\alpha(x,y,n).$

For any two points $x$ and $y$ from $S^N$ we shall denote by
$\gamma(x,y)$ spherical distance between these two points.
Actually, $\gamma(x,y)$ is a measure of angle between $x$ and $y.$
It is obvious, that $\gamma(x,y)\leq \pi.$

We proceed to the formulation of the fundamental results of the
paper.

\begin{Theorem}\label{T1}
Let $\Theta^\alpha(x,y,n)$ is the kernel of Riesz means of the
spectral
expansions\\
 1)if $|\frac{\pi}{2}-\gamma|<\frac{n}{n+1}\frac{\pi}{2}$ then we have
$$
\Theta^\alpha(x,y,n)=O(1)(\frac{n^{(N-1)/2-\alpha}}{(\sin\gamma)^{(N-1)/2}(\sin(\gamma/2))^{1+\alpha}}+
\frac{n^{(N-3)/2-\alpha}}{(\sin\gamma)^{(N+1)/2}(\sin(\gamma/2))^{1+\alpha}}+$$

$$ +\frac{n^{-1}}{(\sin(\gamma)/2)^{1+N}});
$$

2) if $0\leq\gamma\leq\pi$, then we have
$$\Theta^\alpha(x,y,n)=O(1)n^N;$$

 3) if $0<\gamma_0\leq\gamma\leq\pi$, then we have $$\Theta^\alpha(x,y,n)=O(1)n^{N-\alpha};$$
\end{Theorem}

\section{Preliminaries}

In this section we give some properties of Riesz means of spectral
expansions.

Let $f$ be a function with the support in $(0,+\infty).$ If $f$
has the local bounded variation, then the Riesz means $f^\alpha$
for all $Re(\alpha)>-1$ is defined by
\begin{equation}\label{P1}
f^\alpha(t)=\int\limits_0^t\left(1-\frac{s}{t}\right)^\alpha
df(s).
\end{equation}
Using the properties of the distribution $R_\alpha(t)=\alpha
t^{\alpha-1}$, in \cite{HR1} was established the generalization of
M.Riesz theorem :

\begin{Theorem}\label{HR}
Let $\zeta$ is complex number with $Re(\zeta)>0,$ and $M_0(t)$ and
$M_1(t)$ are positive increasing functions in $(0,+\infty)$. If
$f(t)=0, t<0$, and for all $t>0$ satisfied inequalities:

\begin{equation}
|f(t)|\leq M_0(t)
\end{equation}

\begin{equation}
|t^\zeta f^\zeta(t)|\leq M_1(t)
\end{equation}

then for $0<Re(\alpha)<Re(\zeta)$ we have

\begin{equation}
|t^\alpha f^\alpha(t)|\leq
C(1+|\alpha|)^{Re(\zeta)+2}\left(\frac{|\alpha|}{Re(\alpha)}+\frac{|\zeta-\alpha|}
{Re(\zeta-\alpha)}\right)
M_0^{\frac{Re(\zeta-\alpha)}{Re(\zeta)}}(t)M_1^\frac{Re(\alpha)}{Re(\zeta)}.
\end{equation}

Constant $C$ depends only on $\zeta$.

\end{Theorem}

We are going to prove the Theorem \ref{T1} using this statement on
Riesz means.

Let us introduce the Cesaro means of spectral expansions. The
Cesaro means of order $\alpha\geq 0$ is defined as
\begin{equation}\label{CM}
C_n^\alpha
f(x)=\sum\limits_{k=0}^n\frac{A_{n-k}^\alpha}{A_n^\alpha}Y_k(f,x)
\end{equation}
where $Y_k(f,x)\in H_k,k\geq 0$ and
$A_m^\alpha=\frac{\Gamma(\alpha+m+1)}{\Gamma(\alpha+1)m!},m=0,1,2,...,
\Gamma(z)$ is Gamma function.

As a Riesz means the Cesaro means are integral operator with the
kernel
\begin{equation}\label{SC}
\Xi^\alpha(x,y,n)=\sum\limits_{k=0}^n\frac{A_{n-k}^\alpha}{A_n^\alpha}Z_k(x,y),
\end{equation}
where $Z_k$ is zonal harmonics of order $k$, which is reproducing
kernel for the space $H_k,k\geq 0.$ For the kernel we have

\begin{Lemma}\label{KG}
If $\alpha>-1$ and $|\frac{\pi}{2}-\gamma(x,y)|\leq
\frac{n}{n+1}\frac{\pi}{2}$, then
$$
\frac{\Gamma(\alpha+1)}{\Gamma(n+\alpha+1)}\frac{\Gamma(n+\frac{N+1}{2})}{\Gamma(\frac{N+1}{2})}
\frac{\sin\left((n+\frac{N+1}{2})\gamma-(\frac{N-1}{2}+\frac{\alpha}{2})\frac{\pi}{2}\right)}
{(2\sin\gamma)^{(N-1)/2}(2\sin(\gamma/2))^{1+\alpha}}+$$
$$+
\frac{O\left(n^{(N-3)/2}\right)}{(\sin\gamma)^{(N+1)/2}(\sin(\gamma/2))^{1+\alpha}}+
\frac{O\left(n^{-1}\right)}{(\sin(\gamma)/2)^{1+N}});
$$
If $\alpha>-1$ and $0<\gamma_0\leq\gamma\leq\pi,$ then for all
$n>1$
$$
|\Xi^\alpha(x,y,n)|\leq c n^{N-1-\alpha},
$$
If $\alpha>-1$ and $0<\gamma_0\leq\gamma\leq\pi,$ then for all
$n>1$
$$
|\Xi^\alpha(x,y,n)|\leq c n^N.
$$
\end{Lemma}

This lemma is proved in \cite{KG}. It is well known that for
integer order $\alpha$ kernels $\Theta^\alpha(x,y,n)$ and
$\Xi^\alpha(x,y,n)$ have same asymptotical behavior. Our purpose
to extend this result for all $\alpha$.

Let us denote by
$$M_\alpha(t)=\frac{t^{(N-1)/2}}{(\sin\gamma)^{(N-1)/2}(\sin(\gamma/2))^{1+\alpha}}+
\frac{t^{(N-3)/2}}{(\sin\gamma)^{(N+1)/2}(\sin(\gamma/2))^{1+\alpha}}
+\frac{t^{-1}}{(\sin(\gamma)/2)^{1+N}};
$$
If $\zeta $ is integer, then for the kernel $\Theta^\zeta(x,y,n)$
we have
$$
|\Theta^0(x,y,n)|\leq M_0(t);
$$
and

$$
t^\zeta|\Theta^\zeta(x,y,n)|\leq M_\zeta(t).
$$

So using the statement of Theorem \ref{HR} we obtain
\begin{equation}
|t^\alpha \Theta^\alpha(t)|\leq
C(1+|\alpha|)^{Re(\zeta)+2}\left(\frac{|\alpha|}{Re(\alpha)}+\frac{|\zeta-\alpha|}
{Re(\zeta-\alpha)}\right)
M_0^{\frac{Re(\zeta-\alpha)}{Re(\zeta)}}(t)M_1^\frac{Re(\alpha)}{Re(\zeta)}.
\end{equation}
So finally by simplifying we get
$$
t^\alpha\Theta^\alpha(x,y,n)\leq
\frac{t^{(N-1)/2}}{(\sin\gamma)^{(N-1)/2}(\sin(\gamma/2))^{1+\alpha}}+
\frac{t^{(N-3)/2}}{(\sin\gamma)^{(N+1)/2}(\sin(\gamma/2))^{1+\alpha}}+$$

$$ +\frac{t^{-1}}{(\sin(\gamma)/2)^{1+N}};
$$

\section {The estimates for maximal operators}
In this section we are going to prove the estimation for maximal
operators of the Riesz means by using the results of previous
section. Let us recall some standard definitions from harmonic
analysis on unit sphere.

 Spherical ball $B(x,r)$ of radius $r$ and with the
center at point $x$ defined by $B(x,r)=\{y\in S^N:
\gamma(x,y)<r\}.$ For integrable function $f(x)$ the maximal
function of Hardy-Littlwood
\begin{equation}\label{max}
f^*(x)=\sup_{r>0}\frac{1}{|B(x,r)|}\int\limits_{S^N}|f(y)|d\sigma(y)
\end{equation}
is finite almost everywhere on sphere. The maximal function $f^*$
plays a major role in analysis and has been much studied
(see.\cite{ST5}). In particular, for any $p>1$ and if $f\in L_p$,
then there exists constant $c_p$, such that
$$
\|f^*\|_{L_p}\leq\ \frac{c_p(N)}{p-1}\|f\|_{L_p},
$$
where $c_p$ has no singularities at point $p=1.$

\begin{Theorem}\label{T3}
Let $\alpha>\frac{N-1}{2}$ then for all $f\in L_1(S^N)$ we have
\begin{equation}\label{T1}
E_*^\alpha f(x)\leq\frac{c_\alpha(N)}{\alpha-\frac{N-1}{2}}
\left(f^*(x)+f^*(\bar{x})\right)
\end{equation}
\end{Theorem}
where $\gamma(x,\bar{x})=\pi$.

Proof. Since Riesz means of the Fourier-Laplace series are
integral operator with the kernel $\Theta^\alpha(x,y,n)$, we use
the asymptotic behavior of this kernel in Theorem\ref{T1}.

Using this estimates for the kernel $\Theta^\alpha(x,y,n)$ we can
estimate the Riesz means of the spectral expansions:
$$E_n^\alpha f(x)=\int\limits_{S^N}\Theta^\alpha(x,y,n)f(y)d\sigma(y)$$

separate into four part as follow
$$E_n^\alpha f(x)=\int\limits_{S^N}\Theta^\alpha(x,y,n)f(y)d\sigma(y)=$$
$$
\int\limits_{\gamma(x,y)<\frac{1}{n}}\Theta^\alpha(x,y,n)f(y)d\sigma(y)+
+\int\limits_{\frac{1}{n}<\gamma(x,y)\leq\frac{\pi}{2}}\Theta^\alpha(x,y,n)f(y)d\sigma(y)+$$
$$+\int\limits_{\frac{\pi}{2}<\gamma(x,y)\leq\pi-\frac{1}{n}}\Theta^\alpha(x,y,n)f(y)d\sigma(y)+
\int\limits_{\pi-\frac{1}{n}<\gamma(x,y)\leq\pi}\Theta^\alpha(x,y,n)f(y)d\sigma(y).$$
 and for estimate each part let us apply the
Theorem \ref{T1}.

 $$|E_n^\alpha f(x)|\leq C\left(n^N\int\limits_{\gamma(x,y)<\frac{1}{n}}|f(y)|d\sigma(y)+
 n^{\frac{N-1}{2}-\alpha}\int\limits_{\frac{1}{n}<\gamma(x,y)\leq\frac{\pi}{2}}
(\sin\gamma)^{-\frac{N+1}{2}-\alpha}|f(y)|d\sigma(y)\right.+$$

$$+n^{\frac{N-3}{2}-\alpha}\int\limits_{\frac{1}{n}<\gamma(x,y)\leq\frac{\pi}{2}}
(\sin\gamma)^{-\frac{N+3}{2}-\alpha}|f(y)|d\sigma(y)+
n^{-1}\int\limits_{\frac{1}{n}<\gamma(x,y)\leq\frac{\pi}{2}}(\sin\gamma)^{-1-N}|f(y)|d\sigma(y)+$$

$$+n^N\int\limits_{\pi-\frac{1}{n}<\gamma(x,y)\leq\pi}|f(y)|d\sigma(y)+
n^{\frac{N-1}{2}-\alpha}\int\limits_{\pi/2<\gamma(x,y)\leq\pi-\frac{1}{n}}
(\sin\gamma)^{-\frac{N+1}{2}-\alpha}|f(y)|d\sigma(y)+$$

$$\left.+n^{\frac{N-3}{2}-\alpha}\int\limits_{\pi/2<\gamma(xy)\leq\pi-\frac{1}{n}}
(\sin\gamma)^{-\frac{N+3}{2}-\alpha}|f(y)|d\sigma(y)+
+n^{-1}\int\limits_{\pi/2<\gamma(x,y)\leq\pi-\frac{1}{n}}(\sin\gamma)^{-1-N}|f(y)|d\sigma(y)\right).
$$

Let us denote by $U_n(x)$ and $V_n(x)$ the first four member and
the last four member, respectively:

$$
U_n(x)=n^N\int\limits_{\gamma(x,y)<\frac{1}{n}}|f(y)|d\sigma(y)+
n^{\frac{N-1}{2}-\alpha}\int\limits_{\frac{1}{n}<\gamma(x,y)\leq\frac{\pi}{2}}
(\sin\gamma)^{-\frac{N+1}{2}-\alpha}|f(y)|d\sigma(y)+
$$

$$
+n^{\frac{N-3}{2}-\alpha}\int\limits_{\frac{1}{n}<\gamma(x,y)\leq\frac{\pi}{2}}
(\sin\gamma)^{-\frac{N+3}{2}-\alpha}|f(y)|d\sigma(y)+
n^{-1}\int\limits_{\frac{1}{n}<\gamma(x,y)\leq\frac{\pi}{2}}(\sin\gamma)^{-1-N}|f(y)|d\sigma(y)
$$

and

$$
V_n(x)=n^N\int\limits_{\pi-\frac{1}{n}<\gamma(x,y)\leq\pi}|f(y)|d\sigma(y)+
n^{\frac{N-1}{2}-\alpha}\int\limits_{\pi/2<\gamma(x,y)\leq\pi-\frac{1}{n}}
(\sin\gamma)^{-\frac{N+1}{2}-\alpha}|f(y)|d\sigma(y)+
$$

$$
+n^{\frac{N-3}{2}-\alpha}\int\limits_{\pi/2<\gamma(xy)\leq\pi-\frac{1}{n}}
(\sin\gamma)^{-\frac{N+3}{2}-\alpha}|f(y)|d\sigma(y)+
n^{-1}\int\limits_{\pi/2<\gamma(x,y)\leq\pi-\frac{1}{n}}(\sin\gamma)^{-1-N}|f(y)|d\sigma(y)
$$
It is not hard to see that $U_n(\overline{x})=V_n(x)$, where
$\overline{x}$ os opposite point to $x\in S^N,$ i.e.
$\gamma(x,\overline{x})=\pi.$

If we define function $F(t)$ by
$$
F(t)=\int\limits_{\gamma(x,y)<t}|f(y)|d\sigma(y)
$$
then it is easy to show that
$$
F(t)\leq C\ t^Nf^*(x)
$$

Due to definition of $F(t)$ we can rewrite the expression of
$U_n(x)$ as follow:

$$
U_n(x)=n^NF(1/n)+
n^{\frac{N-1}{2}-\alpha}\int\limits_{\frac{1}{n}}^{\frac{\pi}{2}}
(\sin t)^{-\frac{N+1}{2}-\alpha}F'(t)dt+
$$
$$+n^{\frac{N-3}{2}-\alpha}\int\limits_{\frac{1}{n}}^{\frac{\pi}{2}}
(\sin t)^{-\frac{N+3}{2}-\alpha}F'(t)dt+
n^{-1}\int\limits_{\frac{1}{n}}^{\frac{\pi}{2}}(\sin
t)^{-1-N}F'(t)dt.
$$
Integrating by parts, we have

$$
U_n(x)\leq C\ f^*(x)
\left(1+n^{\frac{N-1}{2}-\alpha}\int\limits_{1/n}^{\pi/2}
t^{-\frac{N+3}{2}-\alpha}dt+n^{\frac{N-3}{2}-\alpha}\int\limits_{1/n}^{\pi/2}
t^{-\frac{N+5}{2}-\alpha}dt +n^{-1}\int\limits_{1/n}^{\pi/2}
\frac{dt}{t^2}\right)
$$
The interior of right part of the last inequality may be compute
exactly and

$$U_n(x)\leq \frac{C_1}{\alpha-\frac{N-1}{2}}f^*(x)
$$

Analogously for $V_n(x)$ we have

$$
V_n(x)\leq \frac{C_2}{\alpha-\frac{N-1}{2}}f^*(\overline{x})
$$

 So finally for Riesz means we have

So providing the inequality $\alpha>\frac{N-1}{2}$ we have
\begin{equation}
E_*^\alpha f(x)\leq
\frac{C_\alpha(N)}{\alpha-\frac{N-1}{2}}\left(f^*(x)+f^*(\overline{x})\right).
\end{equation}

Theorem \ref{T3} is proved.

\end{document}